\documentclass[review]{elsarticle}

\usepackage{lineno,hyperref}
\usepackage{amsmath,amsthm,amssymb}
\modulolinenumbers[5]

\journal{Journal of \LaTeX\ Templates}

\newtheorem{thm}{Theorem}
\newtheorem{cl}{Claim}

\newdefinition{rmk}{Remark}
\newproof{pf}{Proof}
\newproof{pot}{Proof of Claim}

\begin{document}

\begin{frontmatter}

\title{Decomposing Claw-free Subcubic Graphs and $4$-Chordal Subcubic Graphs}


\author[1]{Elham Aboomahigir}
\ead{mahigir.elham@gmail.com}

\author[2]{Milad Ahanjideh\corref{mycorrespondingauthor}}
\cortext[mycorrespondingauthor]{Corresponding author}
\ead{ahanjidm@gmail.com}

\author[3]{Saieed Akbari}
\ead{s{\_}akbari@sharif.edu}

\address[1]{Farhangian University, Tehran, Iran}
\address[2]{No 22, 4/1 Alley, East Molavi St., Shahrekord, Chaharmahal and Bakhtiari , Iran}
\address[3]{Department of Mathematical Sciences, Sharif University of Technology, Tehran, Iran}

\begin{abstract}
Hoffmann-Ostenhof's conjecture states that the edge set of every connected cubic graph can be decomposed into a spanning tree, a matching and a $2$-regular subgraph. In this paper, we show that the conjecture holds for claw-free subcubic graphs and $4$-chordal subcubic graphs.
\end{abstract}

\begin{keyword}
\texttt Subcubic graph \sep Hoffmann-Ostenhof’s conjecture \sep Claw-free graph \sep $4$-Chordal subcubic graph\sep Planar graph
\MSC[2010] 05C70
\end{keyword}

\end{frontmatter}

\section{Introduction}\label{sec1}
In this paper, all graphs are assumed to be finite without loops or multiple edges. Let $G$ be a finite  graph with the vertex set $V(G)$ and the edge set $E(G)$. For a vertex $v\in V(G)$, the degree of $v$ in $G$ and the maximum degree of $G$ are denoted by $d_G(v)$ and $\bigtriangleup(G)$, respectively. Here $N_G(v)$ denotes the set of all neighbours of $v$.  The complete graph of order $n$ is denoted by $K_n$.  The complete bipartite graph with partite sets of sizes $m$ and $n$ is denoted by $K_{m, n}$.  A graph is called \textit{cubic} if the degree of every vertex is $3$ and it is called a \textit{subcubic} graph if its maximum degree is at most $3$.
   A graph is called \textit{claw-free} if it has no induced subgraph isomorphic to $K_{1,3}$. A cycle is called \textit{chordless} if it has no chord. A graph $G$ is called \textit{chordal} if every cycle of $G$ of length greater than $3$ has a chord and a graph is \textit{$4$-chordal} if it has no induced cycle of length greater than $4$.  A \textit{cut-edge} of a connected graph $G$ is an edge $e\in E(G)$ such that $G\setminus e$ is disconnected. A subdivision of a graph $G$ is a graph obtained from $G$ by replacing some of the edges of $G$ by internally vertex-disjoint paths. An edge decomposition of a graph $G$ is called a {\it $3$-decomposition}, if  the edges of $G$ can be decomposed into a spanning tree, a matching and a $2$-regular subgraph (the matching or the $2$-regular subgraph may be empty).
   
Hoffmann-Ostenhof proposed the following conjecture in his thesis \cite{Hoff1},  this conjecture also was appeared as a problem of BCC22 \cite{Cameron}.
\vspace{2mm}

\noindent \textbf{Hoffmann-Ostenhof's conjecture.} Every connected cubic graph admits a 3-decomposition.

\vspace{2mm}
 Hoffmann-Ostenhof's conjecture is known to be true for some families of cubic graphs. Kostochka \cite{Kostochka} showed that the Petersen graph, the prism over cycles, and many other graphs have $3$-decompositions. Bachstein \cite{Bachstein} proved that every $3$-connected cubic graph embedded in torus or Klein-bottle has  a $3$-decomposition. Furthermore, Ozeki and Ye \cite{Ozeki} proved that $3$-connected  plane cubic  graphs have $3$-decompositions. Akbari et. al. \cite{Akbari} showed that Hamiltonian cubic graphs have 3-decompositions. Also, it has been proved that the traceable cubic graphs have $3$-decompositions \cite{Abdolhosseini}. In $2017$, Hoffmann-Ostenhof et. al. \cite{Hoff3} proved that planar cubic graphs have $3$-decompositions. In this paper it is shown that the connected claw-free subcubic graphs and the connected $4$-chordal subcubic graphs have $3$-decompositions.

\section{Connected Claw-free Subcubic Graphs Have $3$-Decompositions}
In this section, we show that Hoffmann-Ostenhof's conjecture holds even for  claw-free subcubic graphs.
\begin{thm}\label{A}
If $G$ is a connected claw-free subcubic graph, then $G$ has a $3$-decomposition.
\end{thm} 
\begin{pf}
We apply by induction on $n=|V(G)|$. If $|V(G)|\leq 3$, then the assertion is trivial. Now, we consider the following two cases:
\\
\noindent\textbf{Case 1.} Assume that the graph $G$ has a cut-edge $e$. Indeed we would like to prove that a connected subcubic graph $G$ (not necessarily claw-free) with a cut-edge $e$ has a $3$-decomposition if and only if each component of $G\setminus e$ has a $3$-decomposition.
\noindent
Let $H$ and $K$ be the connected components of $G\setminus e$. By induction hypothesis, both $H$ and $K$ have $3$-decompositions. Let $T_i$, $i=1, 2$, be the spanning trees in the $3$-decompositions of $H$ and $K$. Add $e$ to $T_1\cup T_2$ and consider this tree as the spanning tree in a $3$-decomposition of $G$. Note that we  take the union of cycles and matchings obtained in two components $H$ and $K$ as the $2$-regular subgraph and the matching in the $3$-decomposition of $G$. Now, if $G$ has a $3$-decomposition, then $e$ is contained in the spanning  tree  $T$. Since $T\setminus e$ is union of two trees which
one of them is a spanning tree of $H$ and another is a spanning tree of $K$, so we are done.
\\
\noindent\textbf{Case 2.} By Case 1 we may assume that $G$ is $2$-edge connected. 
If $G$ is triangle-free, then since $G$ is claw-free, we conclude that $\bigtriangleup (G)\leq 2$. Thus $G$ is a cycle and hence in this case the assertion is trivial. \\
\noindent Now, let $xyzx$ be a triangle in $G$. If $d_G(x)=d_G(y)=2$ and $d_G(z)=3$, then $G$ has a cut-edge, a contradiction.\\
 Assume that $d_G(x)=d_G(y)=d_G(z)=3$.
If there is a vertex incident to all $x,y$ and $z$, then $G=K_4$ which satisfied in the conjecture.  Since $G$ is $2$-edge connected, $H=G\setminus \{x, y, z\}$ is connected. By induction hypothesis, $H$ admits a $3$-decomposition. Let $e, f$ and $g$ be the three edges with one end-point in $\{x,y,z\}$, and another end-point in $V(H)$. Add $e$, $f$ and $g$ to $T_1$ to obtain a spanning tree for $G$, where $T_1$ is the spanning tree in the $3$-decomposition of $H$. Now, consider $xyzx$ as a cycle in the $3$-decomposition of $G$, as desired.\\
Finally assume that $d_G(x)=2$  and $d_G(y)=d_G(z)=3$. Assume that $y$ and $z$ have a common neighbour, say $b\neq x$. Clearly, if $d_G(b)=2$, then $G$ has a $3$-decomposition. Now, if $d_G(b)=3$, then $G$ has a cut-edge, a contradiction. Now, assume that $N_G(y)\cap N_G(z)= \{x\}$. Identify $x, y$ and $z$. Call the new vertex $a$ and denote the resulting graph by $H'$. Clearly, $H'$ is a claw-free subcubic graph. By induction hypothesis $H'$ has a $3$-decomposition. Let $T_1$ be the spanning tree in the $3$-decomposition of $H'$. If $d_{T_1}(a)=2$, then let $T$ be the spanning tree of  the $3$-decomposition of $G$ formed by $T_1\cup \{xy, yz\}$. Also consider $xz$ as an edge of the matching in the $3$-decomposition of $G$. Finally, if  $d_{T_1}(a)=1$, then the edge incident with $z$ is contained in $T_1$. Let $T$ be a spanning tree of $G$ formed by $T_1\cup\{xy, yz\}$. Note that $xz$ is contained in the matching in the $3$-decomposition of $G$. The proof is complete.\qed
\end{pf} 
\section{Connected $4$-Chordal Subcubic Graphs Have $3$-Decompositions}
In this section we show that every connected $4$-chordal subcubic graph has a $3$-decomposition.   
\begin{thm}
If $G$ is a connected $4$-chordal subcubic graph, then $G$ has a $3$-decomposition.
\end{thm}
\begin{pf}
 We show that if $G$ does not have a $3$-decomposition, then $G$ is planar and so  by Corollary $13$ of \cite{Hoff3}, $G$  has  a $3$-decomposition, as desired.  To the contrary, 
 suppose that  $G$ is not planar. By Kuratowski's Theorem \cite[p.310]{West}, $G$ either contains a subdivision of $K_5$ or a subdivision of $K_{3,3}$. Since every vertex of $K_5$ has degree $4$, $G$ cannot contain a subdivision of $K_5$. Now, suppose that $G$ contains a subdivision of $K_{3,3}$. \\
 First, we introduce some notation which we need for the rest of the proof.\\
 Let $K_{3,3}^*$ be a subgraph of $G$ which is a subdivision of $K_{3,3}$.
For every edge $e\in E(K_{3,3})$, let $L(e)$ be the set of all new added vertices on the edge $e$ in $G$. For every $v\in V(K_{3,3})$, denote the set of all edges incident with $v$ by $\{e_i(v)\mid 1\leq i\leq 3\}$. If for every $e\in E(K_{3,3})$, $L(e)=\varnothing$, then $G= K_{3,3}$ and so $G$ has a $3$-decomposition, as desired. The rest of the proof is based on the following three claims:
\begin{cl}\label{1}
For every $e\in E(K_{3,3})$, $|L(e)|\geq 2$.
\end{cl}
\begin{pot}
By contradiction, if there exists an edge $e\in E(K_{3,3})$ such that $L(e)=\{ p \}$ and $e$ is incident with the vertex $a$ in $K_{3,3}^*$, then let $C$ be a $4$-cycle of $K_{3,3}$ containing $e$. Let $C^*$ be the subdivision of $C$ in $K_{3,3}^*$. Assume that $\langle C^*\rangle$ denotes the induced subgraph of $G$ on $V(C^*)$. If $\langle C^*\rangle$ is chordless, then $G$ has an induced cycle of length at least $5$, a contradiction. Now, let $wx$ be a chord of $\langle C^*\rangle$ in $G$ such that the cycle $a\ldots wx\ldots upa$ is chordless in $G$ (note that there might exist some vertices between $a$ and $w$ and  also between $x$ and $u$), see Figure 1. So $G$ contains  an induced cycle of length at least $5$,  a contradiction.  If $x=p$, then consider the path between $p$ and $w$ on $C^*$ which contains $a$ and remove all of its vertices except $p$ and $w$. Similarly, one can see that the remaining graph has
an induced cycle of length at least $5$, a contradiction. So the claim is proved. \qed
\end{pot}

\begin{figure}[h!]
\centering
\includegraphics[scale=0.3]{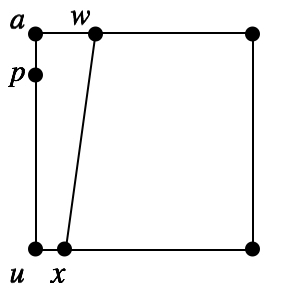}
\caption{The subdivision of  a $4$-cycle in $K_{3,3}^*$}
\end{figure}

For every $v\in V(K_{3,3})$, by Claim 1, let $M(e_{i}(v))$ be the set of two vertices in $L(e_{i}(v))$ which have the smallest distances from $v$. In particular, for $a\in V(K_{3,3})$, let $M(e_1(a))=\{p,q\}$, $M(e_2(a))=\{r,s\}$ and $M(e_3(a))=\{t,z\}$, see Figure 2. 
\begin{figure}[h!]
\centering
\includegraphics[scale=0.3]{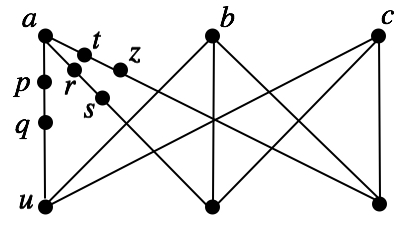}
 \caption{The subdivision of $K_{3,3}$ in $G$}
 \label{22}
\end{figure}

\vspace{2mm}

\begin{cl}\label{2}
 If there exists $v\in \{a,b,c\}$ ($a,b,c$ are the vertices of one part of $K_{3,3}$)  and $y\in M(e_i(v))$,  for some $i$, $1\leq i\leq 3$, such that  $N_G(y)\cap M(e_j(v))=\varnothing$ for every $j\in \{1, 2, 3\}\setminus\{i\}$, $1 \leq j\leq 3 $, then $G$ has a $3$-decomposition.
 \end{cl}
\begin{pot}
Assume that $v=a$ and $y=p$. Since $d_G(q)\leq 3$, there exists $j$, $j\in\{2, 3\}$, such that $N_G(q)\cap M(e_j(a))=\varnothing$. Now, consider a $4$-cycle in $K_{3,3}$ containing $e_1(a)$ and $e_j(a)$ and call it by $C$. Without loss of generality, let $j=3$.
 Assume that $C^*$ is the subdivision of $C$ in $K_{3,3}^*$.  If $\langle C^*\rangle$ is chordless, then by Claim 1,  $C^*$ is an induced cycle of length at least $12$, a contradiction. Now, suppose that the cycle $\langle C^*\rangle$ has at least one chord. Let $wx$ be that chord in $\langle C^*\rangle$ such that the cycle  $atz\ldots wx\ldots qpa$ is chordless (note that $w\in\{t,z\}$ is possible). This implies that $G$ has an induced cycle of length at least $5$, see Figure 3, a contradiction. The same conclusion can be deduced for the case $v=a$ and $y=q$.\qed
 \end{pot}
\begin{figure}[h!]
\centering
\includegraphics[scale=0.3]{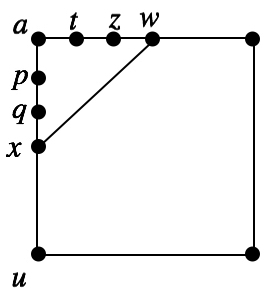}
\caption{The subdivision of  a $4$-cycle in $K_{3,3}^*$}
\end{figure}
\vspace{2mm}
 
\begin{cl}
If  for every $v\in \{a, b, c\}$ ($a,b,c$ are the vertices of one part of $K_{3,3}$) and any $y\in M(e_i(v))$, $i=1,2,3$, $$N_G(y)\cap (\cup_{j\in \{1,2,3\}\setminus \{i\}} M(e_j(v))\neq \varnothing,$$ then $G$  has  a $3$-decomposition. 
\end{cl}
\begin{pot}
Without loss of generality assume that $v=a$ and $N_G(p)\cap M(e_2(a))= \{r\}$. If $N_G(q)\cap M(e_3(a))= \{z\}$, then $G$ contains the induced $5$-cycle $apqzta$, a contradiction. Now,  assume that $N_G(q)\cap M(e_3(a))= \{t\}$. Note that $s$ and $z$ are adjacent. In this case by considering  the induced $5$-cycle $arszta$, we obtain a contradiction to the fact that $G$ is $4$-chordal.\qed
\end{pot}
Thus $G$ is planar and  by Corollary $13$ of \cite{Hoff3}, $G$  has  a $3$-decomposition.\qed
\end{pf}


{\small
\begin{thebibliography}{999}
\bibitem{Abdolhosseini}  F.~Abdolhosseini,  S.~Akbari,  H.~Hashemi,  M.S.~Moradian, Hoffmann-Ostenhof’s conjecture for traceable cubic graphs, (2016), arXiv:1607.04768v1.
\bibitem{Akbari}  S.~Akbari, T.R.~Jensen,  M.~Siggers, Decomposition of graphs into trees, forests, and regular subgraphs, Discrete Math. 338 (8)(2015), 1322--1327. 
\bibitem{Bachstein}   A.C.~Bachstein, Decomposition of Cubic Graphs on The Torus and Klein Bottle, A Thesis Presented to the Faculty of the Department of Mathematical Sciences Middle Tennessee State University, 2015. 
\bibitem{Cameron}  P.J.~Cameron, Research problems from the BCC22, Discrete Math.
311 (2011), 1074--1083.
\bibitem{Hoff1}  A.~Hoffmann-Ostenhof, Nowhere-zero Flows and Structures in Cubic Graphs, Ph.D. thesis (2011), University of Vienna.
\bibitem{Hoff3}  A.~Hoffmann-Ostenhof,  T.~Kaiser,  K. Ozeki, Decomposing planar cubic graphs, 
J. Graph Theory. 88 (2018), 631--640. 
\bibitem{Kostochka}  A.~Kostochka, Spanning trees in $3$-regular graphs, REGS in Combinatorics, University of Illinois at Urbana Champaign, (2009). http://www.math.uiuc.edu/~west/regs/span3reg.html 
52 (2016), 40--46.
\bibitem{Ozeki} K.~Ozeki, D.~Ye, Decomposing plane cubic graphs, European J. Combin.  52 (2016), 40--46. 
\bibitem{West} D.B.~West, Introduction to Graph Theory, Prentice-Hall (2001).




\end{thebibliography}}
\end{document}